\documentclass[letterpaper,12pt,oneside]{article}
\usepackage{amsmath,amsthm}
\usepackage{amsfonts}
\usepackage{latexsym}
\usepackage{xcolor}
\usepackage{multicol}
\usepackage{enumitem}
\usepackage[latin1]{inputenc}
\usepackage{amssymb}
\usepackage{oldgerm}
\usepackage{array}
\usepackage{hyperref}
\theoremstyle{plain}

\newtheorem{lema}{Lemma}
\newtheorem{teor}{Theorem}

\newtheorem{rem}{Remark}

\def\C{\mathbb{C}}

\title{On the separation of the roots of the generalized Fibonacci polynomial}

\author{
Jonathan Garc\'ia \footnote{Departamento de Matem\'aticas, Universidad del Valle, 25360 Cali, Calle 13 No 100-00, Colombia. 
e-mail: garcia.jonathan@correounivalle.edu.co}\,\,\,\,
Carlos A. G\'omez \footnote{Departamento de Matem\'aticas, Universidad del Valle, 25360 Cali, Calle 13 No 100-00, Colombia; Research Group ALTENUA: Algebra, Teoría de Números y Aplicaciones; Colciencias code: COL0017217.
e-mail: carlos.a.gomez@correounivalle.edu.co} \,\,\,\, 
Florian Luca \footnote{School of Mathematics, University of the Witwatersrand, Johannesburg, South Africa; Research Group in Algebraic Structures and Applications, King Abdulaziz University, Jeddah, Saudi Arabia; Max Planck Institute for Software Systems Saarbr\"ucken, Germany and Centro de Ciencias Matem\'aticas UNAM, Morelia, Mexico.  
		e-mail: florian.luca@wits.ac.za}}

\begin{document}

\maketitle

\begin{abstract}
In this paper we prove some separation results for the roots of the generalized Fibonacci polynomials and their absolute values. 

\end{abstract}



\section{Introduction}

A sequence $(u_n)_{n\in {\mathbb Z}} \subseteq \mathbb{C} $ is a linear recurrence sequence of order  $ k\in\mathbb{Z}^{+} $ if it satisfies the recurrence relation 
$$
u_{n+k}=a_{1}u_{n+k-1}+a_{2}u_{n+k-2}+\cdots+a_{k}u_{n} \quad  {\rm for ~all} ~~ n\geq 0
$$ 
with coefficients $ a_1,\ldots,a_k\in \mathbb{C} $ and $a_k\ne 0$. We assume that $k$ is minimal with the above property. Such a sequence $(u_{n})_{n\in {\mathbb Z}} $ has an associated characteristic polynomial given by $f(X) = X^{k} -a_{1}X^{k-1}-\cdots-a_{k}$. Let $\alpha_1, \ldots, \alpha_s \in \C$ be the distinct roots of $\Psi_{k}(X)$. From the theory of linear recurrence sequences (see~\cite[Theorem C.1]{Shorey86}) there are complex single--variable polynomials $A_1, A_2, \ldots, A_s$, uniquely determined by the initial values $u_0, \ldots, u_{k-1}$, such that for all $n\in {\mathbb Z}$
$$
u_n = A_1(n)\alpha_1^{n} + A_2(n)\alpha_2^{n} + \cdots + A_s(n)\alpha_s^n.
$$
From here we see the importance of knowing some properties of the roots of $f(X)$ if one wants to deduce arithmetic properties for the members of $(u_n)_{n\in {\mathbb Z}}$.
 
A special and recently studied case is related to the roots of the characteristic polynomial of the {\it $k$--generalized Fibonacci sequence}\footnote{For an integer $k \geq 2$, the $k$--{\it Fibonacci sequence} is defined by the recurrence relation  $F_n^{(k)}=F_{n-1}^{(k)}+\cdots+F_{n-k}^{(k)}\quad \text{for all}\quad n\geq 2$ and initial values $F_{i}^{(k)}=0$  for $i=-(k-2),\ldots,0$, followed by $F_{1}^{(k)}=1$. The Fibonacci numbers are obtained for $k=2$.}. Let $f_k(X)$ be the generalized Fibonacci polynomial given by
$$
f_k(X)=X^k-X^{k-1}-\cdots-X-1.
$$
This polynomial has exactly one root outside the unit disk. It is real and larger than $1$ and we denote it by $\alpha_k$  (see \cite{DAW}) . 
Let $\rho_1>\cdots>\rho_K$ be the absolute values of the remaining roots of $f_k(X)$ which are smaller than $1$. Then $K=\lfloor (k-1)/2\rfloor$. Let $\alpha_j=\rho_j e^{i\theta_j}$, $\theta_j\in (0,\pi]$ for $j=1,\ldots,K$ be roots of 
$f_k(X)$ which are inside the upper half ${\text{\rm Im}}(z)\ge 0$ of the disk $|z|\le 1$. Then all roots of $f_k(X)$ are 
$$
\rho_1 e^{i\theta_1}, \rho_1e^{-i \theta_1},\ldots,\rho_{K} e^{i \theta_K}, \rho_{K} e^{-i \theta_K},~\alpha_k,\qquad  {\text{\rm when}}\quad k~{\text{\rm is odd}}
$$ 
and 
$$
\rho_1 e^{i\theta_1}, \rho_1e^{-i \theta_1},\ldots,-\rho_K,~\alpha_k,\qquad {\text{\rm when}}\quad k~{\text{\rm is~even}}.
$$
The following bounds on $\rho_j$ for $j=1,\ldots,K$ appear in \cite{CMUC}.

\begin{lema}
	\label{lem:rho}
	We have
	$$
	1-\frac{\log 3}{k}<\rho_j<1-\frac{1}{2^8 k^3}\qquad {\text{\rm for}}\qquad j=1,\ldots,K.
	$$
\end{lema}

The distribution of the arguments $\theta_j$ is also understood. The following appears in \cite{AKLS}, responding to a conjecture proposed in \cite{CMUC}.

\begin{lema}
	\label{lem:theta}
	Let $\alpha_j=\rho_j e^{i \theta_j}$ with $\theta_j\in [0,2\pi)$ for $j=1,\ldots,k$ be all the roots of $f_k(X)$. Then for every $h\in \{0,1,\ldots,k-1\}$, there exists $j$ such that
	$$
	\left|\theta_j-\frac{2\pi h}{k}\right|<\frac{\pi}{k}.
	$$
\end{lema}

To reconcile with our notation, the angles $\theta_j$ for $j=K+1,\ldots,k$ are chosen so that $\theta_{j+\ell}=\theta_j+\pi$ for $\ell=1,\ldots,K$ when $k$ is odd. When $k$ is even, then 
the same formula holds for $\ell=1,\ldots,K-1$ and $\theta_K=\pi$. Finally, $\theta_k=0$ corresponds to the dominant root $\alpha_k$.  

\begin{rem}
	\label{rem1}
	Since the intervals 
	\begin{equation}
		\label{eq:interval}
		\left(\frac{2h \pi}{k}-\frac{\pi}{k},\frac{2h \pi}{k}+\frac{\pi}{k}\right)
	\end{equation}
	are disjoint modulo $2\pi$ as $h$ ranges in $\{0,\ldots,k-1\}$, and since $f_k(X)$ has $k$ roots, each interval above contains the argument of only one of the roots of $f_k(X)$. 
\end{rem}

Another important and useful property is related to the separation of the roots of $f_k(X)$. In the previous paper \cite{GoLu} we proved that
$$
\frac{\rho_i}{\rho_j}>1+8^{-k^4}\qquad {\text{\rm for~all}}\quad 1\le i<j\le K,
$$
and Dubickas \cite{Du}  improved the right--hand side above to  $1+1.454^{-k^3}$.

Our first result is improving this bound.

\begin{teor}
	\label{thm:1}
	The inequality
	\begin{equation}
		\label{eq:lowb}
		\frac{\rho_i}{\rho_j}>1+\frac{1}{10k^{9.6} (\pi/e)^k}\quad {\text{\rm holds~for~all}}\quad 1\le i<j\le K.
	\end{equation}
	and all $k\ge 4$.
\end{teor}

\section{An auxiliar separation result}

The above Theorem 1 is a separation result concerning the absolute values of the differences of the roots of $f_k(X)$. Quite general separation results of this kind appear in \cite{Du} but they are much worse (the denominator of the analogous expression from the right--hand side in \eqref{eq:lowb} in \cite{Du} is exponential in $k^2$). To prove Theorem \ref{thm:1} we need a good separation result on the $\alpha_i$'s themselves.

Before presenting our result on the root separation of $f_k(X)$, with which we will prove Theorem 1, we will show how we can obtain a better result than the one proved by Dubickas \cite{Du} in our particular case by using results of Mahler \cite{Mah} and Mignotte \cite{Mig}. Namely, let us show that the inequality 
\begin{equation}
\label{eq:weak}
|\alpha_i-\alpha_j|>\frac{1}{k^{3/2} 3^{k/2}}\quad {\text{\rm holds~for}}\quad 1\le i<j\le k
\end{equation}
provided $k\ge 100$. This can be proved using an off--the-shelf result. Namely, let 
\begin{equation*}
g_k(X)=(X-1)f_k(X)=X^{k+1}-2X^k+1 \qquad {\rm and} \qquad h_k(X)=X^{k+1}g_k(1/X).
\end{equation*}
Then the roots of $h_k(X)$ are $1$ and $1/\alpha_i$ for $i=1,\ldots,k$. Let them be $y_1,\ldots,y_{k+1}$ with $y_{\ell} =1/\alpha_{\ell}$ and $\alpha_{k+1} = 1$. 
	By \cite{Mah} and \cite{Mig}, the inequality 
	$$
	|y_i-y_j|>\frac{\sqrt{3 |{\text{\rm disc}}(h_k)|}}{d^{d/2+1} \| h_k\|_2^{d-1}}\qquad {\text{\rm holds~for~all}}\quad 1\le i<j\le k+1,
	$$
	holds where $d={\text{\rm deg}}(h_k)$ and ${\text{\rm disc}}(h_k)$ are the degree and discriminant of $h_k(X)$, respectively. 
	For us, $d=k+1$, and 
	\begin{eqnarray}
		\label{eq:calcdisc}
		|{\text{\rm disc}}(h_k)| & = & |{\text{\rm disc}}(g_k)| = (k-1)^2|{\text{\rm disc}}(f_k)| \nonumber\\
      		& = & 2^{k+1} k^k-(k+1)^{k+1}=k^k\left(2^{k+1}-(k+1)\left(1+\frac{1}{k}\right)^k\right)\nonumber\\
		& > & k^k(2^{k+1}-e(k+1))>2^k k^k
	\end{eqnarray}
	since $2^{k}>e(k+1)$ holds for $k\ge 100$. Further, $\|h_k\|_2={\sqrt{1+2^2+1}}={\sqrt{6}}$. We thus get
	\begin{eqnarray*}
		\label{eq:3}
		|y_i-y_j| & > & \frac{{\sqrt{3\cdot 2^k k^k}}}{(k+1)^{(k+3)/2} {\sqrt{6}}^k}=\frac{3^{1/2}}{(k+1)^{3/2}} \frac{1}{(1+1/k)^{k/2}} \left(\frac{1}{3^{k/2}}\right)\nonumber\\
		& \ge & \left(\frac{3}{e}\right)^{1/2}  \frac{1}{(k+1)^{3/2}\cdot 3^{k/2}}>\frac{1.05}{k^{3/2} 3^{k/2}},
	\end{eqnarray*}
	where we used the fact that $(3/e)^{1/2}(1/(k+1)^{3/2})>1.05/k^{3/2}$ for $k\ge 100$. Evaluating the above in $y_i=1/\alpha_i,~y_j=1/\alpha_j$, we get that
	$$
	|\alpha_i-\alpha_j|>\frac{1.05 |\alpha_i\alpha_j|}{k^{3/2} 3^{k/2}}>\frac{1}{k^{3/2} 3^{k/2}},
	$$
	where we used the fact that 
	$$
	1.05|\alpha_i\alpha_j|>1.05\left(1-\frac{\log 3}{k}\right)^2>1
	$$
	(see Lemma \ref{lem:rho}) for $k\ge 100$.
 We next prove a better result.

\begin{teor}
	\label{thm:2}
	The inequality
	$$
	|\alpha_i-\alpha_j|>\frac{1}{k^{6.6} (\pi/e)^{k}}\qquad {\text{\rm holds~for~all}}\quad 1\le i<j\le k
	$$
	and all $k\ge 4$. 
\end{teor}

\section{Proof of the separation result Theorem \ref{thm:2}}

We assume $k\ge 100$. The proof of Mahler's and Mignotte's results from \cite{Mah} and \cite{Mig}, respectively, are based on a discriminant calculation together with an upper bound for a determinant which follows from Hadamard's inequality. Since we know quite a few things about the roots of $f_k(X)$, we visit that proof and use at the appropriate place the extra informations that we have about the roots of $f_k(X)$. We assume that $\alpha_i$ and $\alpha_j$ are 
small roots (i.e., $i,j\in \{1,\ldots,k-1\}$), for otherwise one of $i,j$ is $k$, say $i=k$ and then 
$$
|\alpha_i-\alpha_j|\ge \alpha_k - \rho_j\ge \left(2-\frac{1}{2^{k-1}}\right)-\left(1-\frac{1}{2^8k^3}\right)=1+\frac{1}{2^8k^3}-\frac{1}{2^{k-1}}>0.9,
$$
a much better inequality. We may also assume that $\theta_i\in (0,\pi]$ for otherwise we replace the pair of roots $(\alpha_i,\alpha_j)$ by the pair of roots $({\overline{\alpha_i}},{\overline{\alpha_j}})$ whose separation is the same since $|\alpha_i-\alpha_j|=|{\overline{\alpha_i}}-{\overline{\alpha_j}}|$.  As in the proof of \eqref{eq:weak}, we pass to 
$h_k(X)$ and write 
\begin{equation}
	\label{eq:discform}
	|{\text{\rm disc}}(h_k)|=\prod_{\ell=1}^{k+1} |h_k'(y_\ell)|.
\end{equation}
The left--hand side is 
\begin{equation}
	\label{eq:LHS}
	|{\text{\rm disc}}(h_k)|=2^{k+1} k^k-(k+1)^{k+1}>2^k k^k
\end{equation}
by the calculation \eqref{eq:calcdisc}. In the right--hand side we have
$$
h_k(X)'=(k+1)X^k-2.
$$
When $\ell=k+1$ ($y_\ell=1$), the corresponding factor is 
\begin{equation}
	\label{eq:factor1}
	|h_k'(1)|\le k+3.
\end{equation}
When $\ell=k$ ($y_\ell=1/\alpha_k$), the corresponding factor is   
\begin{equation}
	\label{eq:factor2}
	|h_k'(1/\alpha_k)|=2-\frac{k+1}{\alpha_k^k}<2-\frac{k+1}{2^k}<2\quad {\text{\rm for}}\quad k\ge 100.
\end{equation}
When $y_{\ell}=1/\alpha_{\ell}$ for $\ell\in \{1,\ldots,k-1\}\backslash \{i,j\}$, we have 
\begin{equation}
	\label{eq:factors}
	|h_k'(1/\alpha_{\ell})|\le 2+\frac{k+1}{|\alpha_{\ell}|^k}\le (k+3)|\alpha_{\ell}|^{-k},\qquad \ell\in \{1,\ldots,k-1\}\backslash \{i,j\}.
\end{equation}
Multiplying \eqref{eq:factor1}, \eqref{eq:factor2} and \eqref{eq:factors} for $\ell\in \{1,\ldots,k-1\}\backslash \{i,j\}$, we get
$$
\prod_{\substack{1\le \ell\le k+1\\ \ell\neq i,j}} |h_k'(y_{\ell})|\le 2(k+3)^{k-2} \left|\prod_{\substack{1\le \ell\le k-1\\ \ell\ne i,j}} \alpha_\ell\right|^{-k}.
$$
By the Vieta relations,
$$
\left|\prod_{\substack{1\le \ell\le k-1\\ \ell\ne i,j}} \alpha_\ell\right|^{-1}=\alpha_k |\alpha_i\alpha_j|<\alpha_k<2.
$$
Hence,
\begin{equation}
	\label{eq:RHSpartial}
	\prod_{\substack{1\le \ell\le k+1\\ \ell \neq i,j}} |h_k'(y_{\ell})|<2(k+3)^{k-2} 2^k.
\end{equation}
Now \eqref{eq:discform} together with bounds \eqref{eq:LHS} and \eqref{eq:RHSpartial} give
$$
2^k k^k<|{\text{\rm disc}}(h_k)|=\prod_{\ell=1}^{k+1} |h_k'(y_\ell)|<2|h_k'(y_{i})||h_k'(y_j)| (k+3)^{k-2} 2^k,
$$
which gives 
\begin{equation}
	\label{eq:k/e}
	|h_k'(y_{i})||h_k'(y_j)|>\frac{(k+3)^2}{2(1+3/k)^k}>\frac{(k+3)^2}{2e^3}>\frac{k^2}{2e^3},
\end{equation}
where we used the fact that $(1 +1/x)^x<e$ for $x>1$ with $x=k/3$. 
We work on the left--hand side. We have
$$
|h_k'(y_{i})||h_k'(y_j)|=|y_i-y_j|^2\prod_{\substack{1\le \ell\le k+1\\ \ell\neq i,j}} |y_i-y_{\ell}||y_j-y_{\ell}|.
$$
When $y_{\ell}=1$ or $y_{\ell}=1/\alpha_k$, we have 
$$
|y_i-y_{\ell}||y_j-y_{\ell}|\le (1+1/\rho_i)(1+1/\rho_j)<2.1^2\qquad y_{\ell}\in \{1,1/\alpha_k\}
$$
since $\min\{\rho_i,\rho_j\}\ge 1-\log 3/k>1/1.1$ for $k\ge 100$ by Lemma \eqref{lem:rho}. We thus get that 
\begin{equation}
	\label{eq:RHS2}
	|h_k'(y_{i})||h_k'(y_j)|<|y_i-y_j|^2 (2.1)^4 \prod_{\substack{1\le \ell\le k-1\\ \ell\ne i,j}} |y_i-y_{\ell}||y_j-y_{\ell}|.
\end{equation}
In the right, we write $y_i-y_{\ell}=1/\alpha_i-1/\alpha_{\ell}$ and do the same for $y_j-y_{\ell}$ to get that the right--hand side is 
$$
|\alpha_i-\alpha_j|^2 (2.1)^4 \left(\prod_{\substack{1\le \ell\le k-1\\ \ell\ne i,j}} |\alpha_i-\alpha_{\ell}||\alpha_j-\alpha_{\ell}|\right) |\alpha_i\alpha_j|^{-(k-1)} \prod_{\substack{1\le \ell\le k-1\\ \ell\ne i,j}} |\alpha_{\ell}|^{-2}.
$$
Again by the Vieta relations,
$$
|\alpha_i|^{-2} |\alpha_j|^{-2} \prod_{\substack{1\le \ell\le k-1\\ \ell\ne i,j}} |\alpha_{\ell}|^{-2}=\alpha_k^2<4.
$$ 
Thus, we get that the right--hand side in \eqref{eq:RHS2} is at most 
$$
|\alpha_i-\alpha_j|^2 (2.1)^4\cdot 4 \left(\prod_{\substack{1\le \ell\le k-1\\ \ell\ne i,j}} |\alpha_i-\alpha_{\ell}||\alpha_j-\alpha_{\ell}|\right) |\alpha_i\alpha_j|^{-(k-3)}.
$$
Since 
$$
|\alpha_i|=\rho_i>1-\frac{\log 3}{k}=\exp\left(\log\left(1-\frac{\log 3}{k}\right)\right)>\exp\left(-\frac{2\log 3}{k}\right)
$$
(where we used that $\log(1-x)>-2x$ for $x\in (0,1/2)$), we get that 
$$
|\alpha_i|^{-(k-3)}<|\alpha_i|^{-k}<\exp\left(2\log 3\right)=9,
$$
and the same inequality holds with $i$ replaced by $j$. Thus, the right--hand side of \eqref{eq:RHS2} is at most 
$$
|\alpha_i-\alpha_j|^2 (2.1)^4\cdot 4\cdot 9^2 \left(\prod_{\substack{1\le \ell\le k-1\\ \ell\ne i,j}} |\alpha_i-\alpha_{\ell}||\alpha_j-\alpha_{\ell}|\right).
$$
Since $2.1^4\cdot 4\cdot 9^2<6500$, we get that
$$
|h_k'(y_{i})||h_k'(y_j)|<6500 |\alpha_i-\alpha_j|^2\left(\prod_{\substack{1\le \ell\le k-1\\ \ell\ne i,j}} |\alpha_i-\alpha_{\ell}||\alpha_j-\alpha_{\ell}|\right).
$$
Combining the above with \eqref{eq:k/e}, we get
\begin{equation}
	\label{eq:13000}
	\frac{k^2}{13000 e^3}< |\alpha_i-\alpha_j|^2\left(\prod_{\substack{1\le \ell\le k-1\\ \ell\ne i,j}} |\alpha_i-\alpha_{\ell}||\alpha_j-\alpha_{\ell}|\right).
\end{equation}
It remains to bound the product in the right--hand side.

Let  $h_i,h_j\in \{0,\ldots,k\}$ be such that 
\begin{eqnarray}
\label{cotasthetas}
\theta_i\in \left(\frac{(2h_i-1)\pi}{k},\frac{(2h_i+1)\pi}{k}\right)\quad {\text{\rm and}}\quad \theta_j\in \left(\frac{(2h_j-1)\pi}{k},\frac{(2h_j+1)\pi}{k}\right)
\end{eqnarray}
according to Lemma \ref{lem:theta}. By Remark \ref{rem1}, we have $h_i\ne h_j$. By the same remark, we have that $h_i,h_j$ are both nonzero since the real root $\alpha_k$ corresponds to $h=0$. 
Since $\theta_i\in (0,\pi]$, it follows that $h_i\le (k+1)/2$. We now justify that we may assume that $h_i$ and $h_j$ are consecutive modulo $k$ and since $1\le h_i\le (k+1)/2$, it follows that either $h_j=h_i+1$, or 
$h_j=h_i-1$ and in the second case we must have $h_i\ge 2$. To see why, let us look at $|\theta_i-\theta_j|$. First of all, if $|\theta_i-\theta_j|\ge \pi/2$, it then follows that
\begin{eqnarray*}
|\alpha_i-\alpha_j| & = & \left| \rho_i e^{i\theta_i}-\rho_j e^{i\theta_j}\right| =\left| (e^{i\theta_i}-e^{j\theta_j})+(\rho_i-1)e^{i\theta_i}-(\rho_j-1)e^{i\theta_j}\right| \\
& \ge &  \left| e^{i\theta_i}-e^{j\theta_j}\right| -(1-\rho_i)-(1-\rho_j)\\
& = & 2\left| \sin((\theta_i-\theta_j)/2)\right| - \frac{2\log 3}{k} \ge  {\sqrt{2}} - \frac{2\log 3}{k}>1
\end{eqnarray*}
a much better bound. Thus, we may assume that $\theta_i-\theta_j\in (-\pi/2,\pi/2)$. Now, if $h_i$ and $h_j$ are not consecutive, it then follows that $\theta_i-\theta_j\ge 2\pi/k$. The same calculation as before then gives
\begin{eqnarray*}
|\alpha_i-\alpha_j| &\ge&  2\left| \sin((\theta_i-\theta_j)/2)\right| - \frac{2\log 3}{k}\\
                    &\ge&  (2/\pi) |\theta_i-\theta_j|-\frac{2\log 3}{k}\ge \frac{(4-2\log 3)}{k},
\end{eqnarray*}
again a much better bound. In the above, we used $|\sin x|\ge (2/\pi) |x|$ valid for $x\in (-\pi/2,\pi/2)$ with $x=\theta_i-\theta_j$. 

A bit more can be said. Namely, if $h_j=h_i+1$, then we must have 
\begin{equation}
	\label{eq:hiandhj}
	\frac{(2h_i-0.5)\pi}{k}<\theta_i<\frac{(2h_i+1)\pi}{k}\quad {\text{\rm and}}\quad \frac{(2h_j-1)\pi}{k}<\theta_j<\frac{(2h_j+0.5)\pi}{k}.
\end{equation}
Indeed, say if the first one fails, then $(2h_i-1)\pi/k<\theta_i\le (2h_i-0.5)\pi/k$, while $(2h_j-1)\pi/k=(2h_i+1)\pi/k<\theta_j$, which shows that $|\theta_j-\theta_i|\ge 1.5\pi/k$. As with the previous argument, we get
\begin{eqnarray*}
|\alpha_i-\alpha_j| &\ge& 2\left| \sin((\theta_j-\theta_i)/2) \right| - \frac{2\log 3}{k}\\
                    &\ge& (2/\pi) |\theta_j-\theta_i|-\frac{2\log 3}{k}\ge \frac{3-2\log 3}{k},
\end{eqnarray*}
again a much better inequality. A similar inequality holds if the second inequality  in \eqref{eq:hiandhj} fails. Similar inequalities hold if $h_j=h_i-1$ (just invert the roles of $i$ and $j$ in the above argument to obtain the desired inequalities). From now on we assume that $h_j=h_i+1$ since the other case is obtained by swapping the roles of $i$ and $j$.  Let $\ell\in \{1,\ldots,k-1\}\backslash \{i,j\}$ and let $h_{\ell}$ its corresponding 
$h$ in $\{0,1,\ldots,k-1\}$ according to Lemma \ref{lem:theta}. If $\ell_{i,j}$ is such that $h_{\ell_{i,j}}=h_i-1$, we just bound trivially
\begin{equation}
	\label{eq:triv}
|\alpha_i-\alpha_{\ell}|\le 2.
\end{equation}
If $\ell \neq i,j,\ell_{i,j}$, it then follows that 
\begin{equation}
	\label{eq:w}
	|\theta_i-\theta_{\ell}|\ge \frac{2(|h_\ell-h_i|-1)\pi}{k},
\end{equation}
and $|h_{\ell}-h_i|\ge 2$. 
As $h_{\ell}$ circulates in $\{1,\ldots,k-1\}$ such that $h_{\ell}\not\in \{h_i-1,h_i,h_i+1(=h_j)\}$ the numbers $|h_{\ell}-h_i|-1$ are positive integers $1,2,3,\ldots$ and each one of them is attained at most twice. Assume that 
$\ell$ is such that $\theta_{\ell}$ is ``far" from $\theta_i$, namely $|\theta_i-\theta_{\ell}|\ge \pi/2$. Then 
$$
\left| e^{i\theta_\ell}-e^{i\theta_i} \right| = 2\left| \sin((\theta_{\ell}-\theta_i)/2) \right| \ge {\sqrt{2}},
$$
and 
\begin{eqnarray}
	\label{eq:far}
	|\alpha_i-\alpha_{\ell}| & = & \left| e^{i\theta_i}-e^{i\theta_{\ell}}+(\rho_i-1)e^{i\theta_i}-(\rho_{\ell}-1)e^{i\theta_{\ell}}\right| \le \left| e^{i\theta_i}-e^{i\theta_{\ell}} \right| +\frac{2\log 3}{k}\nonumber\\
	& \le & \left| e^{i\theta_i}-e^{i\theta_{\ell}} \right|  \left(1+\frac{2\log 3}{k |e^{i\theta_i}-e^{i\theta_{\ell}}|}\right)\nonumber\\
	& < & \left| e^{i\theta_i}-e^{i\theta_{\ell}} \right| \left(1+\frac{{\sqrt{2}}\log 3}{k}\right).
\end{eqnarray}
Assume now that $\ell$ is such that $\theta_{\ell}$ is ``close" to $\theta_i$ namely $|\theta_i-\theta_{\ell}|<\pi/2$. Let $w\ge 1$ be such that $|h_{\ell}-h_i|-1=w$. We then have by \eqref{eq:w} that
$$
\frac{\pi}{2}>|\theta_i-\theta_{\ell}|\ge \frac{2w\pi}{k},
$$
so $w \le k/4$. Since every $w$ corresponds to at  most two possible $\ell$'s it follows that there are at most $2\cdot (k/4)=k/2$ such $\ell$'s. For them,
$$
\left| e^{i\theta_i}-e^{i\theta_{\ell}} \right| = 2\left| \sin((\theta_i-\theta_{\ell})/2) \right|  \ge (2/\pi) |\theta_i-\theta_{\ell}|\ge \frac{4 w}{k},
$$
so that 
\begin{eqnarray}
	\label{eq:close}
	|\alpha_i-\alpha_{\ell}| & = & \left| e^{i\theta_i}-e^{i\theta_{\ell}}+(\rho_i-1)e^{i\theta_i}-(\rho_{\ell}-1)e^{i\theta_{\ell}}\right| \le \left| e^{i\theta_i}-e^{i\theta_{\ell}} \right| +\frac{2\log 3}{k}\nonumber\\
	& \le & \left| e^{i\theta_i}-e^{i\theta_{\ell}} \right| \left(1+\frac{2\log 3}{k |e^{i\theta_i}-e^{i\theta_{\ell}}|}\right)\nonumber\\
	& < & \left| e^{i\theta_i}-e^{i\theta_{\ell}}\right| \left(1+\frac{\log 3}{4w}\right).
\end{eqnarray}
As a consequence of \eqref{eq:triv}, \eqref{eq:far} and \eqref{eq:close}, letting $I\le k$ be the number of $\ell$'s for which $h_{\ell}$ is far from $h_i$, we get that
\begin{eqnarray*}
	 \prod_{\substack{1\le \ell\le k-1\\ \ell \ne  i,j}} |\alpha_i-\alpha_{\ell}| &\le& 2\prod_{\substack{1\le \ell\le k-1\\ \ell \ne i,j,\ell_{i,j}}} \left| e^{i\theta_i}-e^{i\theta_{\ell}} \right| \\
	&\times & \left(1+\frac{{\sqrt{2}} \log 3}{k}\right)^I \prod_{1\le w \le k/4} 
	\left(1+\frac{\log 3}{4 w}\right)^2.
\end{eqnarray*}
The first factor in the second line is 
$$
\left(1+\frac{{\sqrt{2}}\log 3}{k}\right)^I<\left(1+\frac{{\sqrt{2}}\log 3}{k}\right)^k<\exp({\sqrt{2}}\log 3)=3^{\sqrt{2}}<5.
$$
The second factor is 
\begin{eqnarray*}
	\prod_{1\le w \le k/4} \left(1+\frac{\log 3}{4 w}\right)^2  < \exp\left(\frac{\log 3}{2} \sum_{w \le k/4} \frac{1}{w}\right) & < & \exp\left(\frac{\log 3}{2}\left(1+\log\left(\frac{k}{4}\right)\right)\right)\\
	& < & k^{\log 3/2}<k^{0.6}.
\end{eqnarray*}
Thus, 
\begin{equation}
\label{eq:10}
\prod_{\substack{1\le \ell\le k-1\\ \ell\ne i,j}} |\alpha_i-\alpha_{\ell}|\le 10 k^{0.6} \prod_{\substack{1\le \ell\le k-1\\ \ell\ne i,j,\ell_{i,j}}} \left| e^{i\theta_i}-e^{i\theta_{\ell}}\right|.
\end{equation}
Finally, 
$$
\left| e^{i\theta_i}-e^{i\theta_{\ell}} \right|  = 2\left| \sin((\theta_i-\theta_{\ell})/2)\right| < |\theta_i-\theta_{\ell}|<\frac{2(|h_i-h_{\ell}|+1)\pi}{k}
$$
where in the last inequality we have used \eqref{cotasthetas} with $\ell$ instead of $j$. As in the previous case, when we were counting $\ell$ such that $\theta_{\ell}$ were close from $\theta_i$, we have that $w:=|h_i-h_{\ell}|$ are now integers of size at most $k/2$, each one of them is counted at twice except that $0$ and $1$ are not counted. 
Thus, we get that 
$$
\prod_{\substack{1\le \ell\le k-1\\ \ell\ne i,j,\ell_{i,j}}} \left| e^{i\theta_i}-e^{i\theta_{\ell}}\right|  \le \left((2\pi/k)^{\lfloor k/2\rfloor}(\lfloor k/2\rfloor+1)!\right)^2 F,
$$
where $F$ accounts for the missing factors (more about that below).  

Now, by Stirling's formula\footnote{$m! =\sqrt{2\pi m}\left(\cfrac{m}{e}\right)^m e^{\lambda_m},\quad {\text{\rm where}}\quad \frac{1}{12m + 1} < \lambda_m <\frac{1}{12m}.
$}: 
\begin{eqnarray*}
(2\pi/k)^{\lfloor k/2\rfloor} (\lfloor k/2\rfloor+1)! &\le& (2\pi/k)^{\lfloor k/2\rfloor} (k/2+1) \lfloor k/2\rfloor !\\
                                                      &\le& 1.1 (k/2+1) (2\pi/k)^{\lfloor k/2\rfloor} \left(\frac{k/2}{e}\right)^{\lfloor k/2\rfloor } {\sqrt{2\pi k/2}}\\
                                                      &\le& 1.1 {\sqrt{\pi}} (k/2+1)  k^{0.5} (\pi/e)^k\\
                                                      & < & k^{1.5} (\pi/e)^{k/2},
\end{eqnarray*}
since $1.1\cdot {\sqrt{\pi}}(k/2+1)<k$ for $k>100$. Thus, 
$$
\prod_{\substack{1\le \ell\le k-1\\ \ell\ne i,j,\ell_{i,j}}} |e^{i\theta_i}-e^{i\theta_{\ell}}|<k^3 (\pi/e)^k F.
$$
On the other hand, $F$ accounts for the fact that when upper bounding the product of $|e^{i\theta_i}-e^{i\theta_{\ell}}|$ by the product of $2w\pi/k$, where $w\le \lfloor k/2\rfloor$ is a positive integer and each $w$ is counted at most twice, we must not count 
$w=1$ and $w=2$ corresponding to $h\in \{h_i-1,h_i,h_i+1\}$. Thus, 
$$
F\le \left(\frac{k}{2\pi}\right)^2 \left(\frac{k}{4\pi}\right)^2=\frac{k^4}{64\pi^4}.
$$
Hence,
$$
\prod_{\substack{1\le \ell\le k-1\\ \ell\ne i,j,\ell_{i,j}}} |e^{i\theta_i}-e^{i\theta_{\ell}}|<k^3 (\pi/e)^k \left(\frac{k^4}{64\pi^4}\right)=\frac{k^7(\pi/e)^k}{64\pi^4},
$$
which together with \eqref{eq:10} gives
\begin{equation*}
	\label{eq:11}
	\prod_{\substack{1\le \ell\le k-1\\ \ell\ne i,j,\ell_{i,j}}} |\alpha_i-\alpha_{\ell}|<\frac{10k^{7.6} (\pi/e)^k}{64\pi^4}.
\end{equation*}
A similar bound applies for $i$ replaced by $j$, so we get that 
$$
\prod_{\substack{1\le \ell\le k-1\\ \ell\ne i,j,\ell_{i,j}}} |\alpha_i-\alpha_{\ell}||\alpha_j-\alpha_{\ell}|<\frac{100 k^{15.2} (\pi/e)^{2k}}{2^{12}\pi^8},
$$
which together with \eqref{eq:13000} gives
$$
\frac{k^2}{13000 e^3}<|\alpha_i-\alpha_j|^2 \left(\frac{100 k^{15.2} (\pi/e)^{2k}}{2^{12}\pi^8}\right),
$$
so
$$
|\alpha_i-\alpha_j|>\frac{1}{k^{6.6} (\pi/e)^k}.
$$
Note that the above inequality was obtained under assumption that $k \le 100$. However, a simple calculation in Mathematica shows that the above inequality holds also for $k<100$, even without the exponential term $\left( \pi/e\right)^k$.  

This completes the proof of Theorem \ref{thm:2}.

\section{The proof of Theorem \ref{thm:1}}

We keep the notations from the previous proof. We may assume that $\theta_i,\theta_j$ are both in $(0,\pi]$, otherwise we replace $\alpha_i$ by ${\overline{\alpha_i}}$ and/or 
$\alpha_j$ by ${\overline{\alpha_j}}$ since this replacement does not change the absolute values $\rho_i$ and $\rho_j$ of the roots. 
Let us assume that the claimed inequality does not hold, namely, that for some $1\le i<j\le K$, we have
$$
\rho_i-\rho_j\le \frac{\rho_j}{10k^{9.6}(\pi/e)^k}<\frac{1}{10k^{9.6}(\pi/e)^k}.
$$
Evaluating $h_k(z)=0$ for $z=y_i$ and its conjugate, we get
$$
y_i^{k+1}=2y_i-1\qquad {\text{\rm and}}\qquad {\overline{y_i}}^{k+1}=2{\overline{y_i}}-1
$$
and multiplying them side by side we get
$$
r_i^{2k+2}=(y_i {\overline{y_i}})^{k+1}=(2y_i-1)(2{\overline{y_i}}-1)=4r_i^2-2(y_i+{\overline{y_i}})+1,
$$
where we write $r_\ell=|y_\ell|=1/\rho_\ell$ for $\ell=1,2,\ldots,K$. A similar relation holds with $i$ replaced by $j$. Subtracting the two relations we get
$$
(r_j-r_i)(r_i^{2k+1}+r_i^{2k}r_j+\cdots+r_j^{2k+1}-4(r_i+r_j))=4({\text{\rm Re}}(y_i)-{\text{\rm Re}}(y_j)).
$$
We discuss the large factor in parenthensis. We have
$$
1<r_{\ell}<\frac{1}{1-\log 3/k}<1+\frac{1.2}{k}\qquad {\text{\rm for}}\quad \ell\in \{i,j\}
$$
by Lemma \ref{lem:rho} since $k>100$. Thus, the factor in parenthesis is at least
$$
2k+2-4(r_i+r_j)>2k+2-4(2+2.4/k)>2k+2-9=2k-7
$$
is positive. On the other hand, this factor in parenthesis is at most as large as
\begin{eqnarray*}
	&& (2k+2)\left(1+\frac{1.2}{k}\right)^{2k+1}\\
	& \le & (2k+2)\left(1+\frac{1.2}{k}\right) \left(\left(1+\frac{1.2}{k}\right)^{k/1.2}\right)^{2/1.2}\\
	& < & 
	\left(2k+2+\frac{2.4k+2.4}{k}\right)e^{2/1.2}<5.3(2k+5)<12k
\end{eqnarray*}
since $k>100$. Further,  
\begin{equation}
\label{eq:rirj}
0<r_j-r_i=\frac{\rho_i-\rho_j}{\rho_i\rho_j}<\frac{r_i}{10 k^{9.6} (\pi/e)^{k}}<\frac{1}{9 k^{9.6}(\pi/e)^k},
\end{equation}
which shows that 
\begin{equation}
\label{eq:Re}
0<{\text{\rm Re}}(y_i)-{\text{\rm Re}}(y_j)<\frac{12 k}{4\cdot 9k^{9.6} (\pi/e)^k}=\frac{1}{3 k^{8.6}(\pi/e)^k}.
\end{equation}
But 
\begin{eqnarray*}
{\text{\rm Re}}(y_i)-{\text{\rm Re}}(y_j) &=& r_i\cos \theta_i -r_j\cos \theta_j\\
                                          &=& r_i\left( \cos \theta_i-\cos \theta_j\right) +(r_j-r_i)\cos\theta_j.
\end{eqnarray*}
Hence, we get from \eqref{eq:rirj} and \eqref{eq:Re} that
\begin{eqnarray*}
	r_i|\cos\theta_i-\cos \theta_j| & < & \left( {\text{\rm Re}}(y_i)-{\text{\rm Re}}(y_j)\right) + (r_j-r_i)\\
	                                & < & \left(3 + \frac{1}{k}\right) \frac{1}{9k^{8.6}(\pi/e)^k}\\
                                 	& < & \frac{1}{2.8k^{8.6}(\pi/e)^k}.
\end{eqnarray*}
Since $r_i>1$, we get that
$$
|\cos\theta_i-\cos \theta_j|<\frac{1}{2.8k^{8.6} (\pi/e)^k}.
$$
We thus get that 
\begin{equation}
	\label{eq:t}
	|\sin((\theta_i-\theta_j)/2)||\sin((\theta_i+\theta_j)/2)|<\frac{1}{2.8k^{8.6} (\pi/e)^k}.
\end{equation}
Fix 
$$
x := (\theta_i+\theta_j)/2, \quad y := (\theta_i-\theta_j)/2 \qquad {\rm so} \qquad \theta_i=x+y,\qquad \theta_j=x-y.
$$ 
Assume that both $\sin(x), \sin(y)$ above are smaller than $1/k^2$. Since $x \in (0,\pi)$, $y \in (-\pi/2,\pi/2)$ and 
$$
\frac{2|y|}{\pi} \le |\sin y|\le \frac{1}{k^2},
$$
we get that $|y|\le ((\pi/2)/k^2)<2/k^2$. Further as in the case of $y$, if $x\in (0,\pi/2)$, then also $x\le 2/k^2$, while if $x\in (\pi/2,\pi)$, we get that $\pi-x<2/k^2$. 

In case $x\in (0,\pi/2)$, so $x\in (0,2/k^2)$, we get that 
$$
\max\{\theta_i, \theta_j\}\le x+|y|\le 4/k^2.
$$ 
Since $\theta_i,\theta_j\in (0,\pi]$, it follows that $\theta_i,\theta_j\in (0,4/k^2)$. In particular, $h_i=h_j=0$, which contradicts Remark \ref{rem1}. 

In the case when $\pi-x<2/k^2$ (i.e. $x \in (\pi/2,\pi)$), we get that $\pi-\theta_{\ell}\in (0,4/k^2)$ holds for both $\ell\in \{i,j\}$. In particular, for $k$ odd the interval corresponding to $h=(k-1)/2$ contains both $\theta_i$ and $\theta_j$, and this is false.    When $k$ is even, we get that the interval corresponding to $h=k/2$ contains both $\theta_i$ and $\theta_j$ and this is again false. 

So, in \eqref{eq:t}, the maximum of the factors is $>1/k^2$, while the minimal factor is $<{\displaystyle{\frac{1}{2.8k^{6.6} (\pi/e)^k}}}$. Assume, for example, that 
$$
|\sin(y)|<\frac{1}{2.8k^{6.6} (\pi/e)^k}.
$$

We now compute, using \eqref{eq:rirj}, 
\begin{eqnarray}
\label{eq:Im}
|{\text{\rm Im}}(y_i)-{\text{\rm Im}}(y_j)| & = & |r_i\sin\theta_i-r_j\sin \theta_j|\nonumber\\
                                            &\le& r_i|\sin \theta_i-\sin\theta_j|+(r_j-r_i)\nonumber\\
                                            & < & 2r_i |\sin(y)\cos(x)|+(r_j-r_i)\nonumber\\
                                         	& \le & \left(2\cdot 1.1 +\frac{2.8}{9k^{3}}\right)\frac{1}{2.8k^{6.6}(\pi/e)^k}\nonumber\\
	& < & \left(\frac{2.3}{2.8}\right)\frac{1}{k^{6.6}(\pi/e)^k}
\end{eqnarray}
for $k>100$. 
It now follows, using \eqref{eq:Re} and \eqref{eq:Im}, that
\begin{eqnarray*}
	|y_i-y_j| & = & {\sqrt{({\text{\rm Re}}(y_i)-{\text{\rm Re}}(y_j))^2+({\text{\rm Im}}(y_i)-{\text{\rm Im}}(y_j))^2}}\nonumber\\
	& < & \left(\frac{1}{9k^4} + \left( \frac{2.3}{2.8}\right) ^2\right)^{1/2}\frac{1}{k^{6.6} (\pi/e)^k}<\frac{1}{k^{6.6}(\pi/e)^k}
\end{eqnarray*}
for $k>100$. In turn, this gives
$$
|\alpha_i-\alpha_j|<\frac{\rho_i\rho_j}{k^{6.6} (\pi/e)^k}
<\frac{1}{k^{6.6}(\pi/e)^k},
$$
contradicting Theorem \ref{thm:2}.  This was when $|\sin((\theta_i-\theta_j)/2)|$ was small. In the case when $|\sin((\theta_i+\theta_j)/2)|$ is small,
the same argument shows that  $|\alpha_i-{\overline{\alpha_j}}|$ is too small again contradicting Theorem \ref{thm:2}. As before, a quick calculation with Mathematica shows that inequality \eqref{eq:lowb} also holds for $k<100$, even without the exponential term $\left( \pi/e\right)^k$.

This finishes the proof of Theorem \ref{thm:1}.

\textbf{Acknowledgement.}
J. G. thanks the Universidad del Valle for support during his master's studies.
C.A.G. was supported in part by Project 71327 (Universidad del Valle).


\begin{thebibliography}{99}

	\bibitem{AKLS} A. Alahmadi, O. Klurman, F. Luca and H. Shoaib, On the arguments of the  roots of the generalized Fibonacci polynomial, {\it Lithuanian Mathematical Journal\/}, to appear.
	
	\bibitem{Du} A. Dubickas, On then distance between two algebraic numbers, {\it Bulletin of the Malaysian Mathematical Science Society\/} {\bf 43} (2020), 3049--3064.
	
	\bibitem{CMUC} C. A. G\'omez and F. Luca, On the distribution of the roots of the polynomial $z^k-z^{k-1}-\cdots-z-1$, {\it Commentaciones Math. Univ. Carolin.\/} {\bf 62} (3) (2021), 291--296.
	
	\bibitem{GoLu} C. A. G\'omez and F. Luca, On the zero--multiplicity of a fifth-order linear recurrence, {\it International Journal of Number
		Theory\/} {\bf  15} (2019), 585--595.
	
	\bibitem{Mah} K. Mahler, An inequality for the discriminant of a polynomial, {\it Mich. Math. J.\/} {\bf 11} (3) (1964), 257.
	
	\bibitem{Mig} M. Mignotte, An inequality about factors of polynomials, {\it Mathematical Computation \/} {\bf 28} (128) (1974), 1153--1157.

\bibitem{Shorey86} T.N. Shorey and R. Tijdeman, {\it Exponential Diophantine equations}, Cambridge University Press, 1986; reprinted 2008.
	
	\bibitem{DAW} D. A. Wolfram, \emph{Solving generalized Fibonacci recurrences}, The Fibonacci Quarterly \textbf{36} (1998), no.~2, 129--145.

\end{thebibliography}
\end{document}